\documentclass[10pt]{amsart}

\usepackage{tikz}
\usetikzlibrary{arrows}

\usepackage{amssymb,amsmath,latexsym,graphicx}
\usepackage{charter,eucal}
\usepackage{amssymb}
\usepackage{amscd}
\usepackage[all,cmtip]{xy}
\usepackage{rotating,multicol}
\usetikzlibrary{decorations.markings}

\newtheorem{theorem}{Theorem }[section]

\newtheorem{lemma}[theorem]{Lemma}
\newtheorem{observation}[theorem]{Observation}

\newtheorem{remark}[theorem]{Remark}
\newtheorem{corollary}[theorem]{Corollary}
\newtheorem{proposition}[theorem]{Proposition}
\newtheorem{principle}[theorem]{\textsc{Principle}}

\newcommand{\bt}{\begin{theorem}}
\newcommand{\et}{\end{theorem}}
\newcommand{\bmt}{\begin{maintheorem}}
\newcommand{\emt}{\end{maintheorem}}
\newcommand{\bc}{\begin{corollary}}
\newcommand{\bl}{\begin{lemma}}
\newcommand{\ec}{\end{corollary}}
\newcommand{\el}{\end{lemma}}
\newcommand{\bo}{\begin{observation}}
\newcommand{\eo}{\end{observation}}
\newcommand{\bp}{\begin{proposition}}
\newcommand{\ep}{\end{proposition}}
\newcommand{\br}{\begin{remark}}
\newcommand{\er}{\end{remark}}
\newcommand{\bpr}{\begin{principle}}
\newcommand{\epr}{\end{principle}}

\def\PG{\mathbf{PG}}

\def\C{\mathbb{C}}
\def\AG{\mathbf{AG}}
\def\AGL{\mathbf{AGL}}

\def\eop{\hspace*{\fill}$\blacksquare$}

\def\bL{\mathbb{L}}

\newcommand{\PGL}{\mathbf{PGL}}
\newcommand{\F}{\mathbb{F}}
\newcommand{\bP}{\mathbb{P}}
\newcommand{\mC}{\mathcal{C}}

\newcommand{\mX}{\mathcal{X}}
\newcommand{\mY}{\mathcal{Y}}

\newcommand{\A}{\mathbb{A}}
\newcommand{\mA}{\mathcal{A}}
\newcommand{\mP}{\mathcal{P}}

\newcommand{\mF}{\mathcal{F}}

\newcommand{\mO}{\mathcal{O}}

\newcommand{\mD}{\mathcal{D}}

\newcommand{\wS}{\texttt{S}}

\newcommand{\Fun}{\mathbb{F}_1}
\newcommand{\Spec}{\texttt{Spec}}
\newcommand{\Proj}{\texttt{Proj}}
\newcommand{\Sch}{\texttt{Sch}}
\newcommand{\Z}{\mathbb{Z}}

\newcommand{\fp}{\frak{p}}
\newcommand{\fq}{\frak{q}}

\title[Projective spaces over $\F_{1^{\ell}}$]{Projective spaces over $\F_{1^{\ell}}$}

\subjclass[2000]{}

\author{Koen Thas}

\thanks{}

\address{Ghent University, Department of Mathematics, Krijgslaan 281, S25, B-9000 Ghent, Belgium}

\email{koen.thas@gmail.com}

\date{}

\begin{document}

\maketitle

\begin{abstract}
In this essay we study various notions of projective space (and other schemes) over $\mathbb{F}_{1^\ell}$, with $\Fun$ denoting the field with one element. Our leading motivation is the ``Hiden Points Principle,'' which shows a huge deviation between the set of rational points as closed points defined over $\mathbb{F}_{1^\ell}$, and the set of rational points defined as morphisms $\Spec(\mathbb{F}_{1^\ell}) \mapsto \mX$.
We also introduce, in the same vein as Kurokawa \cite{Kurozeta}, schemes of $\mathbb{F}_{1^\ell}$-type, and consider their zeta functions. 
\end{abstract}

\setcounter{tocdepth}{1}
\bigskip
{\footnotesize
\tableofcontents
}


\medskip
\section{Introduction}

\subsection{}

The notion of projective and affine geometry plays a crucial role in the emerging theory of $\Fun$, the field with one element. Combinatorially, 
following Tits \cite{anal}, such a projective space $\bP$ of dimension $m$ is a complete graph on $m + 1$ vertices endowed with the full subgraph structure.
Its automorphism group is isomorphic to the symmetric group $\wS_{m + 1}$, and the latter is the Weyl group of the projective general linear group of any projective space of dimension 
$m$ over any field. In \cite{KapranovUN}, Kapranov and Smirnov describe the same theory from the viewpoint of $\Fun$-vector spaces, and they end up with $\wS_{m + 1}$
as the group of $(m + 1)\times(m + 1)$-permutation matrices acting as linear automorphisms.

In Deitmar's scheme theory \cite{Deitmarschemes2}, an affine space (of dimension $m$) over $\Fun$ is defined as $\Spec(\Fun[Y_1,\ldots,Y_m])$,
where $\Fun[Y_1,\ldots,Y_m]$ is the free abelian monoid generated by $Y_1,\ldots,Y_m$ (usually also containing an extra element $0$), and 
the $\Spec$-construction generalizes naturally from commutative unital rings to commutative multiplicative monoids with an extra absorbing element $0$. As noticed by the author in (e.g.) \cite{KT-Chap-Mot} | see also \cite{LPL} |
there is also a natural $\Proj$-construction in this theory, and hence $\Proj(\Fun[Y_1,\ldots,Y_{m + 1}])$ is the ``Algebro-Geom\-etric version'' of $\bP$. 

Both notions fit perfectly together: the $m + 1$ closed points of $\Proj(\Fun[Y_1,\ldots,Y_{m + 1}])$ correspond to the vertices of the graph, and linear subspaces of dimension $1$
indeed contain precisely two closed points.

\subsection{}

In most of the known scheme theories in characteristic $1$, a scheme essentially consists of a triple $(S,X,\gamma)$, where $S$ is a Deitmar scheme, $X$
a Grothendieck scheme, and $\gamma$ a map $S \mapsto X$ which comes from a form of base extension from $\Fun$ to $\Z$. (See for instance Deitmar \cite{Deitmarschemes2}, Connes and Consani \cite{ConCon}, and Thas \cite{KT-Chap-Mot}.) As such, Deitmar schemes play a very important role in Absolute Arithmetic.

In general, many phenomena which are invisible at the $\Fun$-level, only become visible after having applied a base extension functor (say, to a commutative field) | Manin coined this principle with the term ``acquiring flesh.'' This principle is essential for the present note and we study it in the guise of the ``Hidden Points Principle.''

\subsection{}
\label{issue}

In \cite{KapranovUN}, the authors define the ``field extension'' $\F_{1^{\ell}} \Big/ \Fun$ (for any integer $\ell \geq 1$) as the cyclic group $\mu_{\ell}$ with an extra 
element $0$, no addition and usual multiplication. In coordinates, elements of the vector space $V(m + 1,\F_{1^{\ell}})$ are $(m + 1)$-tuples with at most one nonzero entry $a \in \F_{1^{\ell}}$; there is one vector (the zero vector), and
there are $(m + 1)\ell$ vector lines (which are not necessarily different | see \S\S \ref{ghost}), and each contains only the zero vector. Everything passes in a straightforward way to affine spaces, and there arises a multiple occurrence of the aforementioned principle:

\begin{itemize}
\item[{\bf Aff.}]
according to Deitmar's theory, affine $\Fun$-schemes only have one closed point, and the number of lines equals the dimension; 
\item[{\bf Proj.}]
in the $\Proj$-setting, the topology of  $\Proj(\Fun[Y_1,\ldots,Y_{m + 1}])$ and\\ 
$\Proj(\F_{1^{\ell}}[Y_1,\ldots,Y_{m + 1}])$ is the same, so the latter also has 
$m + 1$ closed points.
\end{itemize}

\begin{remark}{\rm
Note that 
\begin{itemize}
\item[$\odot$]
the number of vectors of $V(n,q^\ell)$ with $q \mapsto 1$ is $\lim_{q \mapsto 1}q^{n\ell} = 1$;
\item[$\odot$]
the number of vector lines of $V(n,q^\ell)$ with $q \mapsto 1$ is $\lim_{q \mapsto 1}\frac{q^{n\ell} - 1}{q^\ell - 1} = n$;
\item[$\odot$]
the number of lines of $\AG(n,q^\ell)$ with $q \mapsto 1$ is $\lim_{q \mapsto 1}\frac{q^{n\ell}(q^{n\ell} - 1)}{q^\ell(q^\ell - 1)} = n$.
\end{itemize}
}
\end{remark}

On the other hand, according to Kurokawa \cite{Kurozeta}, a $\Z$-scheme $\mX$ is {\em of $\Fun$-type} if it comes with a counting polynomial $N_{\mX}(T) \in \Z[T]$ (such that for any finite field $\F_q$, we have that $\vert \mX_q \vert = N_{\mX}(q)$). According to Connes and Consani \cite{ConCon} (see also Deitmar \cite{Deitmarschemes2} and  Manin and Marcolli \cite{MM}), we should have that 

\begin{itemize}
\item[{\bf Count.}] 
The equality $\vert \mX_{\F_{1^{\ell}}} \vert = N(\ell + 1)$ holds.
\end{itemize}

So affine and projective spaces over $\F_{1^{\ell}}$ in Deitmar's setting have a number of hidden points which should be manifested through the counting polynomial,
and which are invisible for the Zariski topology at the $\Fun$-level. (Put more bluntly: {\em the combinatorial side and the algebro-geometric side do not agree on field extensions}.)

In this note, we present a way to handle this problem, and show that the algebro-geometric side and the combinatorial one actually {\em do} agree.

\subsection{Acknowledgment}
The author wants to express his gratitude to Nobushige Kurokawa for several very helpful communications about the paper \cite{Kurozeta}.

\medskip
\section{The functor of points}

Let $k$ be a field, and consider an affine algebraic variety $V(\mA)$ with coordinate ring $A = k[t_1,\ldots,t_m]/\mA$. Let $\gamma: k[t_1,\ldots,t_m] \mapsto A$
be the natural projection with kernel $\mA$. Define $u_i := \gamma(t_i)$ for all $i$.
If
\begin{equation}
\varphi:\ A\ \mapsto\ k
\end{equation}
is any morphism, then $(\varphi(u_1),\ldots,\varphi(u_m))$ is a $k$-point of $V(\mA)$. 

Vice versa, let $\overline{x} = (x_1,\ldots,x_m)$ be any $k$-point of $V(\mA)$. Then the following map defines a morphism $A \mapsto k$:
\begin{equation}
\varphi_{\overline{x}}:\ f\ \mapsto\ f(\overline{x}).
\end{equation}

\subsection*{Correspondence}

For any morphism $\varphi: A \mapsto k$, we have that 
\begin{equation}
\varphi \ \equiv \ \varphi_{\overline{x}},
\end{equation}
where $\overline{x} = (\varphi(u_1),\ldots,\varphi(u_m))$.  \\

This classical insight, which sees $k$-points as $k$-morphisms, also works for schemes in general.
This is the motivation for the definition of {\em $\mathbb{F}_{1^\ell}$-rational point} of an $\Fun$-scheme $\mX$, which is used in, e.g., \cite{ConCon}:
\begin{equation}
\mX(\mathbb{F}_{1^\ell}) = \texttt{Hom}(\Spec(\F_{1^\ell}),\mX).
\end{equation}

When $\mX = \Spec(A)$ is an affine $\Fun$-scheme, we also have
\begin{equation}
\mX(\mathbb{F}_{1^\ell}) \cong \texttt{Hom}(A,\F_{1^\ell}).
\end{equation}

\medskip
\section{Spaces over $\F_{1^{\ell}}$}
\label{spa-ext}

It is our goal to find a common picture to the different models of affine and projective $\F_{1^{\ell}}$-spaces. Our starting point is the projection map
\begin{equation}
V(m + 1,\F_{1^{\ell}})\ \ \mapsto\ \ \Proj(\mathbb{F}_{1^\ell}[Y_1,\ldots,Y_{m + 1}]).
\end{equation}

Following Kapranov and Smirnov \cite{KapranovUN} and \cite{KT-Chap-Comb}, $V = V(m + 1,\F_{1})$ has only one vector (corresponding to the zero vector),
and $m + 1$ directions, which after base extension to a field, become vector lines. Each is given a coordinate of type $(-,1,-)$, and the linear maps act on these directions. Defining 
$\PG(m,\Fun)$ from $V$ in the usual way, each direction defines a projective point and the zero vector gets lost, so that automorphisms of $V$ carry over in a faithful manner
to automorphisms of $\PG(m,\Fun)$.

\subsection{Over $\mathbb{F}_{1^{\ell}}$ (ghost directions)}
\label{ghost}

Kapranov and Smirnov consider all directions of the form $(-,\alpha,-)$ with $\alpha \in \mu_{\ell}$ to define $V(m + 1,\F_{1^\ell})$. In fact, as any direction corresponds to a vector line, no new directions are introduced, but rather cyclotomic coordinates get into the picture. Making the transition to projective space as above, we get the same set of  points (see also Proposition \ref{topell1}), as proportionality kills the cyclotomy (on the geometric level). 

If $X$ is a $k$-scheme, where $k$ is a field, then the closed points of $X$ represent orbits of the Galois group $\texttt{Gal}(\overline{k}/k)$, and all $\overline{k}$-rational points of $X \times_{k}\overline{k}$ are contained in the union of these orbits. So there is a natural map
\begin{equation}
\alpha: X(\overline{k}) \longrightarrow X,
\end{equation}
sending closed points of $X(\overline{k})$ to closed points of $X$, which is neither injective nor surjective in general.

Going back to the spaces $\PG(n,\mathbb{F}_{1^\ell})$ of above, we see the closed points of $\PG(n,\mathbb{F}_{1^\ell})$  as orbits of $\texttt{Gal}(\overline{\F_1}/\F_{1^\ell})$.
A stalk at an arbitrary closed point of $\PG(n,\mathbb{F}_{1^\ell})$ is isomorphic to
\begin{equation}
\mu_\ell \times \Fun[X_1,\ldots,X_n].
\end{equation}

So on the algebraic level we can see
the extension of the ground field | we consider the stalk as consisting, besides $0$, of $\ell$ distinct copies of $\Fun[X_1,\ldots,X_n] \setminus \{0\}$ equipped with 
a sharply transitive $\mu_\ell$-action, which is in accordance with the classical picture. Similar conclusions can be made for the affine spaces $\Spec(\F_{1^\ell}[X_1,\ldots,\\ X_n])$ (where the cyclotomic directions are ghost directions defined on the associated vector space, but are not defined in the Deitmar scheme).

For any finite field $\F_p$, where $p = 1$ is allowed, we have a natural identification
\begin{equation}
V(mn,\F_p)\ \overset{\sim}{\mapsto}\ V(m,\F_{p^n}),
\end{equation}
and so we have
\begin{equation}
V(m,\F_{1^\ell})\ \overset{\sim}{\mapsto}\ V(m\ell,\Fun)\ \overset{\otimes \F_p}{\mapsto}\ V(m\ell,\F_p)\ \overset{\sim}{\mapsto}\ V(m,\F_{p^\ell}),
\end{equation}
which justifies Kapranov and Smirnov's approach. Here, the base extension ``$\otimes_{\Fun}\F_p$'' can be defined using, e.g., the affine space: with $A = \F_1[X_1,\ldots,X_{m\ell}]$ we have
\begin{equation}
A \otimes_{\Fun}\F_p \ =\ (A \otimes_{\Fun} \Z) \otimes_{\Z} \F_p\ =\  \Z[A]  \otimes_{\Z} \F_p,
\end{equation}
where $\Z[A]$ is the monoidal ring defined by $A$.

\subsection{Hidden Points Principle}

The following trivial observation reveals the Hidden Points Principle at the topological level.

\begin{proposition}[\cite{KT-Chap-Mot}]
\label{topell1}
Topologically, the structure of $\PG(n,\mathbb{F}_{1^\ell})$ is independent of the choice of $\ell$. 
\end{proposition}

{\em Proof}.\quad
It suffices to observe that for any $\gamma \in \mu_\ell$ and any $i \in \{0,1,\ldots,n\}$, we have
\begin{equation}
(X_i) = (\gamma X_i),
\end{equation}
whence the topology is the same as (= homeomorphic to) that of $\PG(n,\Fun)$.
\eop \\

Due to this fact, we need to introduce ``hidden points,'' which cannot be detected by the classical way in the Zariski topology.


We can formulate the Hidden Points Principle as follows (for a Deitmar scheme $\mX$ defined over $\mathbb{F}_{1^\ell}$): 
\begin{equation}
\vert \texttt{Hom}(\Spec(\F_{1^\ell}),\mX) \vert\ \gg \ \vert \mX_{\mathbb{F}_{1^\ell}} \vert,
\end{equation}
where on the right hand side, $\mX_{\mathbb{F}_{1^\ell}}$ is the set of $\F_{1^\ell}$-rational points defined as closed points (i.e., coming from (next-to-) maximal ideals). In fact, the ratio
\begin{equation}
f(\mX,\ell) := \frac{\vert \texttt{Hom}(\Spec(\F_{1^\ell}),\mX)\vert}{\vert \mX_{\F_{1^\ell}}\vert} 
\end{equation}
is strictly increasing for a fixed $\mX$, as a function of $\ell \in \mathbb{N} \setminus \{0\}$. \\

In this note, we will look for a minimal enlargement of Deitmar schemes which reveals the hidden points in  Deitmar schemes such as  $\Spec(\F_{1^\ell}[X_1,\ldots,X_m])$ and 
$\Proj(\F_{1^\ell}[X_1,\ldots,X_m])$, and which resolves the issues of \S\S \ref{issue}.

\medskip
\subsection{The combinatorial picture | I}

\subsubsection{}

Let $\Gamma = (V,E)$ be a graph (undirected), with $V$ the vertex set and $E$ the edge set. A {\em cover} of $\Gamma$ is a graph $\widetilde{\Gamma} = (\widetilde{V},\widetilde{E})$ together with a surjective graph morphism
\begin{equation}
\gamma: \widetilde{\Gamma} \ \mapsto\ {\Gamma}
\end{equation}
such that locally $\gamma$ is a bijection (that is, for any vertex $v \in \widetilde{V}$, $\gamma$ induces a bijection between the edges incident with $v$ and 
the edges incident with $\gamma(v)$). 

If there is a positive integer $m$ such that any vertex of $\Gamma$ has a fiber of size $m$, then $(\widetilde{\Gamma},\gamma)$ is an {\em $m$-fold cover}.

Now let $A = (V,E)$ and $B = (V',E')$ be two graphs. The {\em lexicographic product} $A[B]$ of $A$ and $B$, also denoted by $A \cdot B$, is the graph defined on the vertex set $V \times V'$, such that $(v,v') \sim (w,w')$ if either $v \sim w$, or ($v = w$ and $v' \sim w'$).  We allow $B$ to be a directed.

The lexicographic product $A \cdot B$ is in general not commutative, and, considered as the morphism defined by
\begin{equation}
\pi: A\ \mapsto\ B:\ (v,w) \ \mapsto\ v,
\end{equation}
not a cover of $A$. On the other hand, each vertex fiber has a constant size $\vert V' \vert$, and the morphism {\em is} surjective. 

Let $G$ be any group, and $S$ a (minimal) generating set of $G$. The {\em Cayley graph} $\Gamma(G,S)$ of $G$ with respect to $S$ is the directed graph
with vertex set $G$, and directed edges of type $(v,sv)$, with $s \in S$. (With any fixed $s$ corresponds a color, but that is not important here.)

\subsection*{Example}

Let $G = \mu_n$ be any cyclic group, with $n \in \mathbb{N}^{\times} \cup \{ \vert \mathbb{N} \vert \}$, and let $S = \{ s\}$, with $s$ a generator. 
Then $\Gamma(\mu_n,S)$ is a directed cycle on $n$ vertices.

\subsubsection{Combinatorial picture}

One way of introducing a combinatorial model for $\bP^m(\F_{1^\ell})$ could be to define it as $K_{m + 1}[\Gamma(\mu_{\ell}),\{s\}]$. This model essentially corresponds to the Kapranov-Smirnov model of vector spaces over extensions of $\Fun$.

\includegraphics{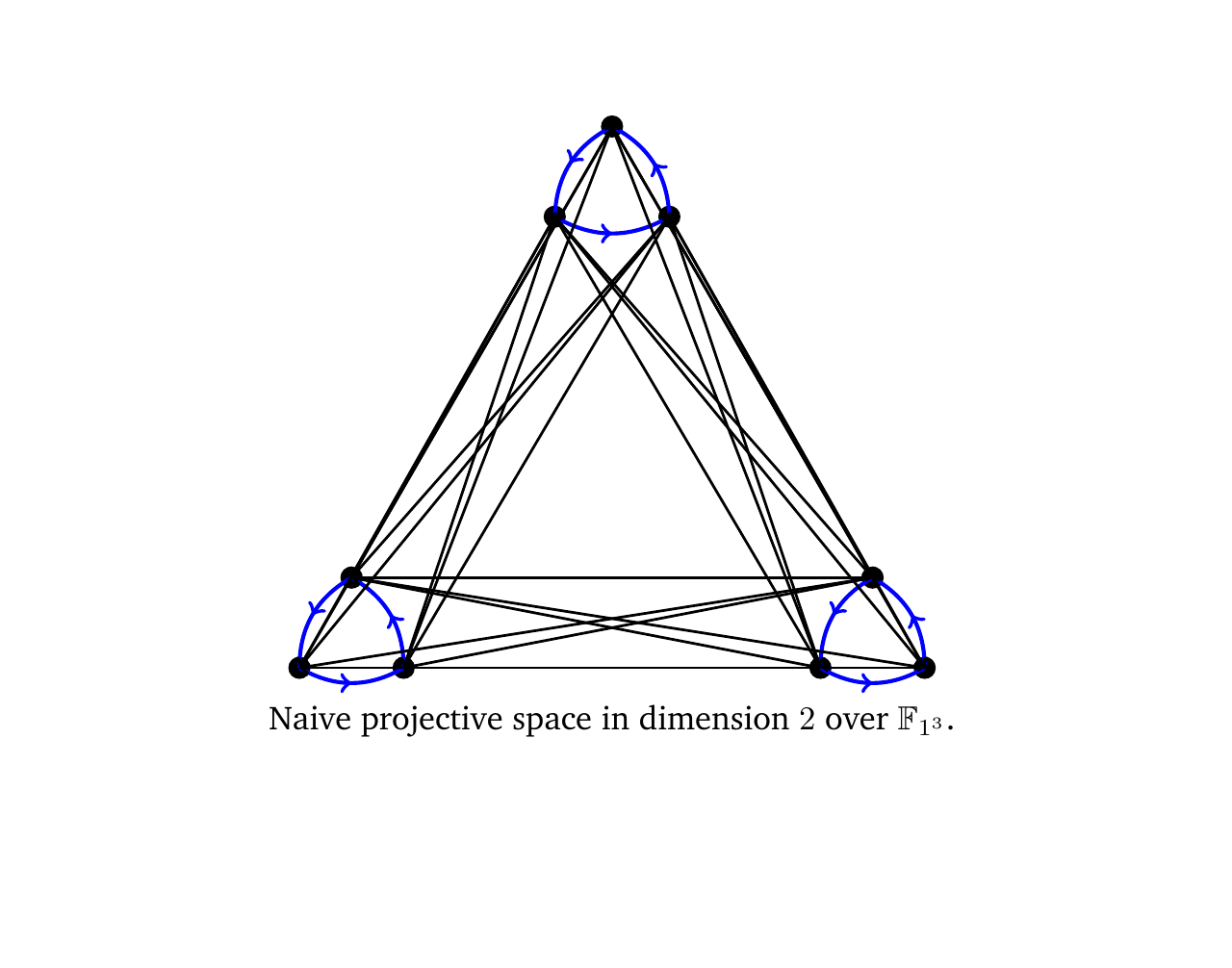}

\subsubsection{Automorphism group}

The automorphism group of 
$K_{m + 1}[\Gamma(\mu_{\ell}),\{s\}]$ is easily seen to be isomorphic to $\mu_{\ell} \wr \wS_{m + 1}$.

The picture presented in this section is certainly richer than the topological model: it now possesses  the right automorphism group, but the counting polynomial is not the predicted one.

\medskip
\section{Frames}
\label{frame}

Let $S$ be any monoid, written multiplicatively and with $0$, and $m \in \mathbb{N}^{\times}$. The {\em affine frame} $\mA(m,S)$ is defined as  the set 
$\{ (s_1,\ldots,s_m) \vert s_i \in S \}$. Let $(0,\ldots,0) =: \omega$.

If $S \setminus \{0\}$ is also a group, the {\em projective frame} $\mP(m - 1,S)$ is defined as $\mA(m,S) \setminus \omega$ modulo proportionality (so an element of $\mP(m - 1,S)$ is  a class $(s_1:\cdots:s_m) := \{ \sigma(s_1,\ldots,s_m) \vert \sigma \in S, \sigma \ne 0  \}$, with $(s_1,\ldots,s_m) \ne \omega$). If $m = 0$, $\mP(m - 1,S) = \emptyset$; if $m = 1$,  $\mP(m - 1,S)$ is a point, and 
\begin{equation}
\mP(1,S) = \{(1 : 0),(0 : 1)\} \cup \{(1 : s) \vert s \in S^{\times}\}. 
\end{equation}

\subsection{Frames over $\F_{1^\ell}$}

Let $\ell \in \mathbb{N}^{\times}$. The affine, respectively projective, frame {\em over} $\F_{1^\ell}$ of {\em dimension} $d$ is $\mA(d,\F_{1^\ell})$, respectively $\mP(d,\F_{1^\ell})$. (Here, $d \in \mathbb{N}$, respectively $d \in \mathbb{N} \cup \{-1\}$.)

We will see $\mA(d,\F_{1^\ell})$ and $\mP(d,\F_{1^{\ell}})$ as models of respectively affine and projective spaces over $\F_{1^\ell}$.

\subsection{Subframes}

An (affine) {\em subframe} of $\mA(d,\F_{1^{\ell}})$ is any subset of tuples $(x_1,\ldots,x_d)$
in which a fixed number $e$ of  coordinates ($0 \leq e \leq d$), are fixed elements of $\F_{1^\ell}$, and the others take all possible values. Its {\em dimension} is $d - e$.

Similarly as in the previous subsection, one defines projective subframes of dimension $d - e - 1$ of $\mP(d - 1,\F_{1^\ell})$ from affine subframes of dimension $d - e$ of $\mA(d,\F_{1^\ell})$ endowed with the  proportionality relation.

\subsection{The combinatorial picture | II. From spaces to frames}

Let $\A^m_{1^\ell} =\\  \Spec(\F_{1^\ell}[X_1,\ldots,X_m])$ be the Deitmar scheme of $m$-dimensional affine space over $\F_{1^\ell}$. According to the point of view of morphisms, a
(rational) point should correspond to a morphism
\begin{equation}
\gamma: \F_{1^\ell}[X_1,\ldots,X_m]\ \mapsto\ \F_{1^\ell}.
\end{equation}

Obviously, any such $\gamma$ is completely determined by $(\gamma(X_1),\ldots,\gamma(X_m))$, and conversely, any such choice gives us  a morphism.
So $\texttt{Hom}(\Spec(\F_{1^\ell}),\A_1^m)$ is in a natural bijective correspondence with the affine frame $\mA(m,\F_{1^\ell})$, and 
\begin{equation}
\vert \texttt{Hom}(\Spec(\F_{1^\ell}),\A_1^m)\vert = (\ell + 1)^m.
\end{equation}

Similar remarks hold for projective space Deitmar schemes and projective frames. In the latter case, we have
\begin{equation}
\vert \texttt{Hom}(\Spec(\F_{1^\ell}),\mathbb{P}_1^m)\vert = \sum_{i = 0}^m(\ell + 1)^i.
\end{equation}

\subsubsection{Example (case $\ell = 1$)}

The number of rational points of $\mP(d,\Fun)$ is $2^{d} + 2^{d - 1} + \cdots + 1$. Amongst these are the $d + 1$ closed points of the underlying Deitmar scheme.

\subsection{Automorphism group}

Obviously, since we do not have any addition in $\F_{1^\ell}$, we have that the automorphism group (defined in a natural way) of $\mA(d,\F_{1^\ell})$, respectively $\mP(d,\F_{1^\ell})$, is
\begin{equation}
\mu_\ell \wr \wS_d, \ \mathrm{respectively}\ \Big(\mu_\ell \wr \wS_{d + 1}\Big) \Big/ \F_{1^\ell}^\times.
\end{equation}

This picture is very interesting: in case of affine and projective spaces over real fields $k$, the respective automorphism groups act transitively on the $k$-linear subspaces of fixed dimension, and in particular on the $k$-points. This is not at all the case here; observe the following. (Denote the affine group by $\AGL_d(1^\ell)$ and the projective group by $\PGL_{d + 1}(1^\ell)$.)

\bo
Let $\mA^r(d,\F_{1^\ell})$ be the set of elements in $\mA(d,\F_{1^\ell})$ which have precisely $r$ nonzero entries ($r \in \{ 0,1,\ldots,d\}$). Similarly we define $\mP^r(d,\F_{1^\ell})$ ($r \in \{1,\ldots,d + 1\}$). Then $\AGL_d(1^\ell)$ stabilizes each set 
$\mA^r(d,\F_{1^\ell})$, and acts transitively on its elements. Also, $\PGL_{d + 1}(1^\ell)$ stabilizes each set 
$\mP^r(d,\F_{1^\ell})$, and acts transitively on its elements. \eop 
\eo

So in this setting, the original closed Deitmar points form one orbit under the full automorphism group. 

\medskip
\br{\rm
If $\mu_\ell \cong \F_q^\times$ for some prime power $q$, then note that the sets $\mA(d,\F_{1^\ell})$ and $\mP(d,\F_{1^\ell})$
coincide with the rational point sets of $\A^d_{\F_q}$ and $\bP^d_{\F_q}$. 
}
\er

\medskip
\section{Frames and enlarged Deitmar schemes}

As we have seen, the data $(K_{m}[\Gamma(\mu_{\ell},S)],\mP(m - 1,\F_{1^\ell}))$ reveals a number of hidden properties of (say) projective spaces $\Proj(\F_{1^\ell}[X_1,\ldots,X_m])$. As the points of an $\F_{1^\ell}$-frame naturally correspond to $\F_{1^\ell}$-rational points of the corresponding space, we want to enlarge Deitmar schemes such that points of the frames associated to affine/projective/etc. spaces appear as closed points.
Those can be seen if we (minimally) enlarge the category of Deitmar schemes in the following way, by using congruences.

\subsection{Congruences and congruence ideals}

Let $M$ be a multiplicative commutative monoid. A {\em congruence} $\mC \subseteq M \times M$ on $M$  is an equivalence relation with the additional property that if $(u,v) \in \mC$, then $(mu,mv) \in \mC$ for each $m \in M$. If $(u,v) \in \mC$, we also write $u \sim v$.
If $\mC$ is a congruence on $M$, one shows that $M/\mC$ (or $M/\sim$) naturally inherits the structure of a monoid. 

Note that in general the union of two congruences $\mC, \mC'$ is not a congruence (as transitivity does not necessary hold anymore). (The congruence generated by $\mC, \mC'$ is essentially the transitive closure.)  The intersection {\em is} a congruence, though.

Any ideal $I$ in the monoid $M$ gives rise to a congruence, called ``Reese congruence" and denoted by $\mC_I$; if $m \in M\setminus I$, the class $[m]$ is just $\{m\}$; if $a, b \in I$, then $[a] = [b]$. So $I$ defines one class, which we will also denote by $[I]$.
If $M$ has an absorbing element $0$, then $0 \in I$, so 
in that case $a \sim 0 \sim b$. From now on we keep assuming that $M$ has a $0$.

A congruence $\mC$ is {\em prime} if from $(ab,ac) \in \mC$ follows that either $(b,c) \in \mC$ or $(a,0) \in \mC$. (Alternatively, a congruence $\mC$ on $M$ is {\em prime} if $M/\mC$ is integral, i.e., if it satisfies the cancellation property.) 
Now let $\fp$ be a prime ideal in the monoid $M$, so that $0 \in \fp$. Let $a, b, c$ be elements in $M$, and suppose $ab = ac$ in $M/\mC_{\fp}$. If $ab \sim ac \in \fp$, then either $a \in \fp$, or $a \not\in \fp$ and $b, c \in \fp$. (If $ab = ac \not\in \fp$, then $a, b, c \ne 0$.) If $M$ is supposed to be integral, then $\mC_\fp$ is a prime ideal of $M$.

\subsection{Example: affine space}

Consider $\Spec(\F_{1^\ell}[X_1,\ldots,X_m])$, Deitmar version. Put $A := \F_{1^\ell}[X_1,\ldots,X_m]$. 
If $\fp$ is a prime ideal, then  $\mC_{\fp}$ is a prime congruence since $A$ is integral. 



\subsection{New Zariski topology $\Spec^c(M)$}

In Deitmar scheme theory, any set of elements $S$ in $\Fun[X_1,\ldots,X_m]$ defines an ideal $I(S)$, and the closed set corresponding to $I(S)$ classically contains the solutions of all the equations defined by its elements. Only equations of the form $P = 0$ arise, with $P$ an element in $\Fun[X_1,\ldots,X_m]$. We need to allow equations of the form $P = 1$, and more generally, over $\F_{1^\ell}$, of the form $P = \mu$.


Let $\Spec^c(M)$ be the set of prime congruences on $M$. Then the Zariski topology it comes with (and denoted in the same way) is defined by the closed sets 
\begin{equation}
C(I) := \{ \fp\ \vert\ I \subseteq \fp, \fp\ \text{prime congruence} \},
\end{equation}
with $I$ any congruence on $M$. 

The key is that the closed set topology is generated by the sets
\begin{equation}
\mD(a,b) := \{ \fp\ \vert\ (a,b) \in \fp, \fp\ \text{prime congruence} \}.
\end{equation}

\subsection{Localization}

Let $\fp$ be any prime congruence. Define $S_{\fp}$ as the monoid multiplicatively generated by the set 
\begin{equation}
\{ b - a\ \vert\ (a,b) \in A \times A,\ a \not\sim_{\fp} b \}. 
\end{equation}

Then $A_{\fp} := S_{\fp}^{-1}A$ is defined as the monoid naturally generated by the fractions $\{ \frac{u}{v}\ \vert\ u \in A, v \in S_{\fp} \}$.
Here, we identify $u/v$ with $u'/v'$ if there is a $w \in S_{\fp}$ such that $wuv' = wu'v$.

\subsection{Sections and structure sheaf}

Let $U$ be an open set. Then $\mO_X(U)$ is defined as the set of maps
\begin{equation}
s : U \ \mapsto\ \coprod_{\fp \in U}A_{\fp}
\end{equation}
such that for each prime $\fp \in U$, $s(\fp) \in A_{\fp}$, and for which there is a neighborhood $V$ of $\fp$ in $U$ such that for 
all $\fq \in V$, $s(\fq)$ is locally a fraction. This means that there are $a \in A$ and $f \in \bigcap_{\fq \in V} S_{\fq}$ so that $s = a/f$ over $V$. 

Note that $\mO_X(U)$ carries the structure of a monoid.

\subsection{Example: maximal ideals (closed points) for affine spaces}

Consider again $A := \F_{1^\ell}[X_1,\ldots,X_m]$. For each coordinate monomial $P$ and each $\mu \in \F_{1^\ell}$, we consider the {\em base congrence} on $A$ generated by
$P \sim \mu$. If, for example, $P = X_j$, the congruence is 
\begin{equation}
\{ (QX_j^n,Q\mu^n) \vert Q \in A, n \in \mathbb{N}\} \cup \{  (Q\mu^n,QX_j^n) \vert Q \in A, n \in \mathbb{N} \}.
\end{equation}  

 It is easy to see that the maximal congruences are precisely those generated by base congruences of the following form
 \begin{equation}
 \{ (X_i,\mu_i)\ \vert\ i \in \{1,\ldots,n \}, \mu_i \in \F_{1^\ell} \}.
 \end{equation}

In other words, the maximal congruences correspond to the rational points/frame points of $\mA(m,\F_{1^\ell})$ as described in \S \ref{frame}.

\subsection{General congruence schemes and $\Proj^c$-schemes}

For the theory of general congruence schemes and more details, we refer to Deitmar's paper \cite{congruence}.

As for the $\Proj^c$-construction, this is tersely described in \cite{MMKT2} in the case of projective spaces. One defines the {\em irreducible congruence} $\texttt{Irr}^c$ on $A = \F_{1^\ell}[X_1,\ldots,X_m]$ as 
\begin{equation}
\texttt{Irr}^c := \Big\langle X_1 \sim 0, \ldots, X_m \sim 0 \Big\rangle,
\end{equation}
leaves it out, and applies a construction which is similar to that in Deitmar scheme theory. The next-to-maximal congruences in $A$ correspond to the closed points, and they correspond to the elements of the projective frame $\mP(m - 1,\F_{1^\ell})$.

\medskip
\subsection{Maximal points vs. closed points}

In the congruence setting of this section, one observes that for affine and projective spaces the set of rational points {\em coincides} with the set  of closed points, a phenomenon which for proper fields only occurs for algebraically closed fields.

\medskip
\subsection{Absolute and geometric $\Fun$-Frobenius endomorphisms}

Consider the algebraic closure $\overline{\mathbb{F}_1}$ of $\F_1$; by definition, it consists of all complex roots of unity, together with an absorbing element $0$, foreseen with multiplication. We define the {\em $\Fun$-Frobenius endomorphism}\index{absolute!Frobenius endomorphism} of degree $n \in \mathbb{N}$, denoted $\mathrm{Fr}^n_1$\index{$\mathrm{Fr}^n_1$}, to be the map
\begin{equation}
\mathrm{Fr}_1^n: \overline{\F_1} \longrightarrow \overline{\F_1}: x \longrightarrow x^n.
\end{equation}

Elements of $\F_{1^{d}} \cong \mu_d \cup \{0\} \leq \overline{\F_1}$ are characterized by the fact that they are the solutions of 
\begin{equation}
{\mathrm{Fr}_1^{d + 1}}(x) = x,
\end{equation}
which is analogous to the fact that elements of finite fields $\F_{q^d} \leq \overline{\F_q}$ are singled out as fixed points of $\mathrm{Fr}^d$.

The absolute $\Fun$-Frobenius map $\mathrm{Fr}_1^{d + 1}$ with $d \in \mathbb{N}^\times$ acts on an extended Deitmar scheme $\mY$ of finite type over $\F_{1^d}$  by acting trivially on the topology, and as $\mathrm{Fr}_1^{d + 1}$ on the  structure scheaf. Passing to the scheme $\overline{\mY} := \mY \otimes_{\Spec(\F_{1^d})}\Spec(\overline{\Fun})$, we define the {\em geometric $\Fun$-Frobenius map} (still in the same degree) as the map 
$\mathrm{Fr}_{\mathrm{geom}}^{d + 1}: \overline{\mY} \mapsto \overline{\mY}$ with local comorphisms $y \otimes f \mapsto y^{d + 1} \otimes f$. It acts 
as 
\begin{equation}
\mathbf{x} \mapsto \mathbf{x}^{d + 1}
\end{equation}
on frame points.

\medskip
\section{Zeta functions}

In the inspiring note \cite{Kurozeta}, Kurokawa says that a $\Z$-scheme $\mX$ is of {\em $\Fun$-type} if there exists a polynomial $N_\mX$ in $\Z[X]$ such that, for each prime power $q$, we have
\begin{equation}
\vert \mX_q \vert = N_{\mX}(q).
\end{equation}
(Through one of Tate's conjectures, such schemes would come with a mixed Tate motive, a feature which one sees more clearly in the appropriate Grothendieck ring | see for instance \cite{KT-Chap-Mot} for a detailed discussion.)

In fact, originally, in \cite{Kurozeta}, a $\Z$-scheme (of finite type) is of $\Fun$-type
if its arithmetic zeta function $\zeta_{\mX}(s)$ can be expressed in the form
\begin{equation}
\zeta_{\mX}(s) = \prod_{k = 0}^n\zeta(s - k)^{a_k}
\end{equation}
with the $a_k$s in $\Z$. (Here, $s$ is a complex variable | cf. the remark right after the statement of Theorem \ref{Kuro}.)
The definition through the counting polynomial then is derived from the following result.

(In the next theorem, $\zeta(\cdot)$ is the classical Riemann zeta.)

\begin{theorem}[Kurokawa \cite{Kurozeta}]
\label{Kuro}
Let $\mX$ be a $\Z$-scheme of finite type. The following are equivalent.
\begin{itemize}
\item[{\rm (i)}]
The arithmetic zeta function of $\mX$ has the following form:
\begin{equation}
\zeta_{\mX}(s) = \prod_{k = 0}^n\zeta(s - k)^{a_k}
\end{equation}
with the $a_k$s in $\Z$.
\item[{\rm (ii)}]
For all primes $p$ we have
\begin{equation}
\zeta_{\mX\vert \F_p}(s) = \prod_{k = 0}^n(1 - p^{k - s})^{-a_k}
\end{equation}
with the $a_k$s in $\Z$.
\item[{\rm (iii)}]
There exists a polynomial $N_{\mX}(Y) = \sum_{k = 0}^na_kY^k$ such that
\begin{equation}
\vert \mX_{\F_{p^m}} \vert = N_{\mX}(p^m) 
\end{equation}
for all finite fields $\F_{p^m}$.
\end{itemize}
\end{theorem}

For details about matters of convergence, we refer to Theorem \ref{eqext} below.

Kurokawa defines the {\em $\Fun$-zeta function}\index{$\Fun$-zeta function} of a $\Z$-scheme $\mX$ of $\Fun$-type as 
\begin{equation}
\zeta_{\mX\vert \Fun}(s) :=  \prod_{k = 0}^n(s - k)^{-a_k}
\end{equation}
with the $a_k$s as above. The {\em Euler characteristic}\index{Euler characteristic} is
\begin{equation}
\vert X_{\Fun} \vert := \sum_{k = 0}^na_k = N_{\mX}(1).
\end{equation}

To fit this formula into {\bf Count}, put 
\begin{equation}
\F_{1^0} = \F_{1^1}.
\end{equation}

(Although this seems rather artificial, it is also rather interesting that one makes {\bf Count} work by putting $0 = 1$!)

\begin{theorem}[Kurokawa \cite{Kurozeta}]
\label{Kurolimit}
Let $\mX$ be a $\Z$-scheme of $\Fun$-type. Then
\begin{equation}
\zeta_{\mX\vert \Fun}(s)\index{$\zeta_{X\vert \Fun}(s)$} =  \lim_{p \longrightarrow 1}\zeta_{\mX\vert \F_p}(s)(p - 1)^{\vert \mX_{\Fun} \vert}.
\end{equation}
Here, $p$ is seen as a complex variable (so that the left hand term is the leading coefficient of the Laurent expansion of $\zeta_{\mX \vert \Fun}(s)$ around $p = 1$).
\end{theorem}

\medskip
\subsection{Schemes of $\Fun$-type vs. schemes defined over $\Fun$}
\label{VS}

As Deitmar envisioned in \cite{Deitmarschemes2}, commutative monoids with a $1$ and $0 \ne 1$ should/could be seen as ``$\Fun$-algebras.'' For each such 
monoid $A$, we have an embedding $\Fun \hookrightarrow A$, so a projection $\Spec(A) \twoheadrightarrow \Spec(\Fun)$. So $\Spec(A)$ is ``defined over $\Fun$.''
(And one can make similar observations for general monoidal schemes.) Considering a $k$-algebra $\mA$ for a field $k$, one can realize $\mA$ as 
a quotient of a polynomial ring over $k$ with the kernel of an appropriate surjective morphism to $\mA$. In particular, if 
\begin{equation}
\iota: k\ \hookrightarrow\ R 
\end{equation}
is an embedding of a field in a unital commutative ring, then $R$ comes with a $k$-algebra structure, so $R$ is a quotient of a polynomial ring over $k$. And 
any unital commutative ring is a $\Z$-algebra, of course.  Over $\Fun$, that is to say, in the setting of monoidal schemes,
this situation is less natural, as some unital commutative monoids  with $0$ can only be 
realized as a quotient of a polynomial ring over $\Fun$ and a {\em congruence} (instead of a monoidal ideal). For instance, although we have an embedding $\Fun \hookrightarrow \Z$, the latter, as a monoid, {\em needs} a congruence which is {\em not} associated to a monoidal ideal if one wants to see it as such a quotient. Whereas over $\F_{1^2}$, $\Z$ is indeed even {\em isomorphic} to a polynomial ring. So it seems more natural to say that $\Spec(\Z)$ is ``defined over $\F_{1^2}$,'' instead than over $\Fun$. But still, strictly speaking, if $\Spec(A)$ is defined over 
$\F_{1^\ell}$, it should also be defined over $\Fun$ as well. With this discussion in mind, one can then decide when a $\Z$-scheme is defined over $\Fun$, or an extension. 

One should make a distinction in any case with Kurokawa's $\Z$-schemes ``of $\Fun$-type.'' If  $\mX$ is a $\Z$-scheme of $\Fun$-type \`{a} la Kurokawa, then any
reasonable definition in the same vein of $\Z$-schemes over extensions of $\Fun$ will lead to the property that ``of $\F_{1^m}$-type'' implies
``of $\F_{1^n}$-type'' if $m$ divides $n$, but not necessarily the other way around. (This phenomenon will manifest itself in the next subsection.)
So the situation is quite different (and more general) than the $\Fun$-algebra viewpoint.

\medskip
\subsection{Dirichlet series}

A {\em Dirichlet series} is a series of the form 
\begin{equation}
D(s) = \sum_{m = 1}^\infty\gamma(m)m^{-s}
\end{equation}
with $s$ a complex variable and $\gamma: \mathbb{N} \to \C$ a function. It can be shown that if $\sum_{m = 1}^\infty\vert \gamma(s)m^{-s}\vert$ does not diverge for all $s$, there exists a real number $\kappa$ such that $\sum_{m = 1}^\infty\vert \gamma(s)m^{-s}\vert$ converges for all $s$ with $\Re(s) > \kappa$ (that is, the Dirichlet series is absolutely convergent for these values).
(The real part of a complex number $c$ is denoted by $\Re(c)$.)

We will use the next well known result later on.

\bt[Uniqueness of coefficients]
Suppose 
\begin{align}
D(s) = \sum_{m = 1}^\infty\gamma(m)m^{-s} \ \text{and}\  E(s) = \sum_{m = 1}^\infty\beta(m)m^{-s}
\end{align}
are Dirichlet series, both absolutely convergent for $\Re(s) > \sigma$. If $D(s) = E(s)$ for each $s$ in some 
infinite sequence $\{s_k\}_k$ such that $\Re(s_k) \to \infty$ as $k \to \infty$, then $\gamma(n) = \beta(n)$ for each $n \in \mathbb{N}$.
\et

\medskip
\subsection{Schemes of $\F_{1^\ell}$-type}

Let $V$ be nonsingular projective variety over $\F_q$, the finite field with $q$ elements. The {\em local zeta function} (or {\em congruence zeta function})
of $V$ is given by
\begin{equation}
\zeta_V(s) = Z(V,s)\ := \ \exp\Big(\sum_{m = 1}^\infty \frac{\vert V_{q^m} \vert}{m}q^{-sm}\Big).
\end{equation}
Sometimes $Z(V,s)$ is also denoted by $Z(V,q^{-s})$.
The variable transformation $u = q^{-s}$ gives $Z(V,u)$ as a formal power series in $u$. To stress that the local zeta 
function is expressed ``with respect to $\F_q$,'' we also write $\zeta_{V \vert \F_q}(s)$ instead of $\zeta_{V}(s)$.

Let $\ell$ be a positive nonzero integer.
We will say that a $\Z$-scheme $\mX$ is {\em of $\F_{1^\ell}$-type} if there exists a polynomial $N_{\mX}(Y) = \sum_{k = 0}^na_kY^k$ such that 
for all primes $p$ and all extensions $\F_{p^{\ell m}}$ of $\F_{p^\ell}$, we have
\begin{equation}
\vert \mX_{\F_{p^{\ell m}}} \vert = N_{\mX}(p^{\ell m}). 
\end{equation}

\bp
If the $\Z$-scheme $\mX$ is of $\F_{1^m}$-type, $m \in \mathbb{N}_0$, then it is also of $\F_{1^r}$-type for any positive integer
multiple $r$ of $m$.
\eop 
\ep

In the next theorem, for a complex number $c$,  the modulus will be denoted by $\vert c\vert$. Also,
$\log(\cdot)$ will denote the principal value complex logarithm.

\begin{theorem}
\label{eqext}
Let $\mX$ be a $\Z$-scheme of finite type, and $\ell$ a positive integer. The following are equivalent.
\begin{itemize}
\item[{\rm (i)}]
There exist $k \in \mathbb{N}$ and integers $a_1,\ldots,a_k$, so that the following equality holds:
\begin{equation}
\prod_{\text{$p$ prime}}\zeta_{\mX\vert\F_{p^\ell}}(s)  = \prod_{k = 0}^n\zeta(\ell(s - k))^{a_k},
\end{equation}
where $s \in \C$ and $\Re(s) > k + \frac{1}{\ell}$.
\item[{\rm (ii)}]
There exist $k \in \mathbb{N}$ and integers $a_1,\ldots,a_k$, so that the following equality holds
for all primes $p$:
\begin{equation}
\zeta_{\mX\vert \F_{p^\ell}}(s) = \prod_{k = 0}^n(1 - p^{\ell(k - s)})^{-a_k},
\end{equation}
where $s \in \C$ and $\Re(s) > k + \frac{1}{\ell}$.
\item[{\rm (iii)}]
There exists a polynomial $N_{\mX}(Y) = \sum_{k = 0}^na_kY^k \in \Z[Y]$ such that
\begin{equation}
\vert \mX_{\F_{p^{\ell m}}} \vert = N_{\mX}(p^{\ell m}) 
\end{equation}
for all finite fields $\F_{p^{\ell m}}$.
\end{itemize}
\end{theorem}

{\em Proof}.\quad
We first show that (ii) and (iii) are equivalent. Observe that, with $q = p^\ell$, 
\begin{equation}
\label{first}
\log(\zeta_{\mX\vert \F_{q}}(s)) = \ \sum_{m = 1}^\infty \frac{\vert \mX_{q^m} \vert}{m}q^{-sm}
\end{equation}
while
\begin{align}
\label{second}
\log\Big(\prod_{k = 0}^n(1 - p^{\ell(k - s)})^{-a_k}\Big) &= \sum_{k = 0}^n-a_k\log(1 - p^{\ell(k - s)}) \nonumber \\
&= \sum_{k = 0}^na_k\Big(\sum_{m = 1}^{\infty}\frac{p^{\ell(k - s)m}}{m} \Big) \nonumber \\
&= \sum_{m = 1}^\infty \frac{\Big(\sum_{k = 0}^na_k{p^{\ell mk}}\Big)}{m}q^{-sm}. 
\end{align}
(Note that we can indeed expand $\log(1 - p^{\ell(k - s)})$ in its Mercator series since $\vert p^{\ell(k - s)} \vert < 1$.)

Then assuming (ii) (equality of the expressions in (\ref{first}) and (\ref{second})), equality of the coefficients in both sides leads to (iii). 
Also, after assuming (iii), substituting $\vert \mX_{\F_{p^{\ell m}}} \vert = N_{\mX}(p^{\ell m})$ for all $p^{\ell m}$ in (\ref{first}) leads to (ii) (using (\ref{second})). 

Now assume (ii). Then
\begin{align}
\prod_{\text{$p$ prime}}\zeta_{\mX\vert\F_{p^\ell}} &= \prod_{\text{$p$ prime}}\Big(\prod_{k = 0}^n(1 - p^{\ell(k - s)})^{-a_k}\Big) \nonumber \\
&= \prod_{k = 0}^n\Big(\prod_{\text{$p$ prime}}(1 - p^{\ell(k - s)})^{-a_k}\Big) \nonumber \\
&= \prod_{k = 0}^n\zeta(\ell(s - k))^{a_k},
\end{align}
which yields (i).\\

Then assume (i). We obtain
\begin{align}
\prod_{\text{$p$ prime}}\zeta_{\mX\vert \F_p^\ell}    &=   \prod_{k = 0}^n\zeta(\ell(s - k))^{a_k} \nonumber \\
&= \prod_{k = 0}^n\Big(\prod_{\text{$p$ prime}}(1 - p^{\ell(k - s)})^{-a_k}\Big) \nonumber \\
&= \prod_{\text{$p$ prime}}\Big(\prod_{k = 0}^n(1 - p^{\ell(k - s)})^{-a_k}\Big). 
\end{align}

Taking complex logarithms of both sides yields: 
\begin{align}
 \sum_{\text{$p$ prime}}\Big(\sum_{m = 1}^\infty \frac{\vert \mX_{q^m} \vert}{m}q^{-sm}\Big) &=
 \sum_{\text{$p$ prime}}\Big(\sum_{k = 0}^n-a_k\log(1 - p^{\ell(k - s)})\Big) \nonumber \\
 &=  \sum_{\text{$p$ prime}}\bigg(\sum_{k = 0}^na_k\Big(\sum_{m = 1}^{\infty}\frac{p^{\ell(k - s)m}}{m} \Big)\bigg). 
 \end{align}
 
 It follows that 
 \begin{align}
 \sum_{p^m\ \text{prime power}}\frac{\vert \mX_{p^{\ell m}} \vert}{m}p^{-\ell ms} &= \sum_{p^m\ \text{prime power}}\Big(\sum_{k = 1}^n a_kp^{\ell mk} \Big)p^{-\ell ms}.
 \end{align}
 
 The uniqueness of the coefficients in the Dirichlet series lead to (iii).
\eop \\

\br{\rm
The reader observes that the switchings above of summation symbols resp. product symbols can indeed be done as $\Re(s) > k + \frac{1}{\ell}$. (For the summation symbols, one can use Fubini's Theorem; switching of the product symbols reduces to summation through application of $\log(\cdot)$.)
}\er

\br{\rm
The author is not aware of a connection between the property expressed in (iii) of the previous theorem and (variations of) mixed Tate motives (as is the case for $\ell = 1$).
}
\er

\br{\rm
Note that the equivalence between (ii) and (iii) remains true if stated for particular primes, and hence particular sets of primes. We will come back to this property further on. 
}\er

Define the {\em $\F_{1^\ell}$-zeta function} of a $\Z$-scheme $\mX$ of $\F_{1^\ell}$-type as 
\begin{equation}
\zeta_{\mX\vert \F_{1^\ell}}(s) :=  \prod_{k = 0}^n(\ell(s - k))^{-a_k}
\end{equation}
with the $a_k$s as in the previous  theorem. 

\medskip
Similarly as Theorem \ref{Kurolimit}, we have:

\begin{theorem}
Let $\mX$ be a $\Z$-scheme of $\F_{1^\ell}$-type. Then
\begin{equation}
\zeta_{\mX\vert \F_{1^\ell}}(s)\index{$\zeta_{X\vert \Fun}(s)$} =  \lim_{p \longrightarrow 1}\zeta_{\mX\vert \F_{p^\ell}}(s)(p - 1)^{\vert \mX_{\Fun} \vert}.
\end{equation}
Here, $p$ is again seen as a complex variable.
\end{theorem}

{\em Proof}.\quad
We have 
\begin{equation}
\zeta_{X\vert \F_{1^\ell}}(s)(p - 1)^{\# X(\Fun)} = \prod_{k = 0}^n\big(\frac{1 - p^{\ell(k - s)}}{p - 1}\big)^{-a_k}
\end{equation}
so that for the limit we get
\begin{eqnarray}
\lim_{p \longrightarrow 1}\zeta_{X\vert \F_{1^\ell}}(s)(p - 1)^{\# X(\Fun)} &= &\prod_{k = 0}^n\big(\lim_{p \longrightarrow 1}\frac{1 - p^{\ell(k - s)}}{p - 1}\big)^{-a_k}\nonumber \\
&= &\prod_{k = 0}^n(\ell(s - k))^{-a_k}\nonumber \\ 
&= &\zeta_{X\vert \F_{1^\ell}}(s). 
\end{eqnarray}
\eop 

\medskip
\section{Schemes of $\Fun$-type versus schemes of $\F_{1^\ell}$-type}

In general, $\Z$-schemes of $\F_{1^\ell}$-type are not necessarily of $\F_{1^m}$-type, with $m$ a natural divisor of $\ell$ (cf. Remark \ref{VS}). In particular, they need not be of $\Fun$-type, as the following example shows. 

\subsection{}

Let $X^2 + uX + v \in \Z[X,X_2,\ldots,X_m]$ and $u$ odd, $\mathrm{gcd}(u,v) = 1$. It is easy to see that over $\F_{p^2}$, $p$ any prime, the number of points is $2(p^2 + 1)^m$. Taking $(u,v)$ such that $X^2 + uX + v$ is not reducible over every prime field $\F_p$, we obtain a scheme of $\F_{1^2}$-type which is not of $\F_1$-type. Its $\F_{1^2}$-counting polynomial is 
\begin{equation}
2(Y + 1)^m = 2\sum_{k = 0}^m\begin{pmatrix}m \\ k \end{pmatrix}Y^{m - k}, 
\end{equation}
so its $\F_{1^2}$-zeta function is
\begin{equation}
\prod_{k = 0}^m(2(s - k))^{-2\begin{pmatrix}m \\ k \end{pmatrix}}.
\end{equation}

\subsection{}

Let $\mX$ be a $\Z$-scheme of finite type defined over $\Fun$. By definition, it means that $\mX$ is a $\Z$-scheme of finite type which arises from Deitmar base extension to $\Z$ from a Deitmar scheme.
Then Deitmar showed in \cite{Deitmarzeta} that there exists a natural number $e$ and a polynomial $N(Y) \in \Z[Y]$ such that for every prime power $q$ for which $\mathrm{gcd}(q - 1,e) = 1$, one has $\#\mX(\F_q) = N(q)$. It might be interesting to consider $\Z$-schemes $\mX$ which come with a counting polynomial $N_{\mX}(Y)$ such that for  some positive integer $\ell$, 
we have that
\begin{equation}
\#\mX(\F_{p^\ell}) = N_{\mX}(p^\ell),
\end{equation}
where $p$ varies over the primes, except possibly a {\em finite number of exceptions}. (When $\ell = 1$, such schemes are related to mixed Tate motives.) \\

Constructing examples is easy.
Consider, for example, the affine conics with equation 
\begin{equation}
\mC(u,v): X^2 - uY^2 = v, 
\end{equation}
with $u, v \in \Z$ (and $u, v$ fixed). Over a finite field $\F_{2^n}$ we have that $\mC(u,v)$ has $2^n + 1$ points. If $p$ is an odd prime, we have that $\vert \mC(u,v)\vert_{\F_{p^n}}$ is $p^n - 1$ if $u$ is a square, and $p^n + 1$ if $u$ is not. If $u$ would be a square in $\F_p$, it is also a square in $\F_{p^2}$, so 
$\vert \mC(u,v) \vert_{\F_p} = p - 1$ and $\vert \mC(u,v) \vert_{\F_{p^2}} = p^2 - 1$. If $u$ is not a square in $\F_p$, 
we have $\vert \mC(u,v)\vert_{\F_p} = p + 1$, but in the quadratic extension $\F_{p^2}$ we know that $u$ {\em is} a square, so that $\vert \mC(u,v) \vert_{\F_{p^2}} = p^2 - 1$.

In other words, we have a fixed counting polynomial for all fields $\F_{p^2}$ ($p$ any prime) in odd characteristic, but not necessarily over $\F_p$. 

\subsection{}

The {\em Grothendieck ring of schemes} of finite type over a field $k$, denoted as $K_0(\Sch_{k})$, is generated by the isomorphism classes of schemes $\mX$ of finite type over $\Fun$, $[\mX]_{k}$, with the relation
\begin{equation}
[\mX]_{k}= [\mX\setminus \mY]_{k} + [\mY]_{k} 
\end{equation}
for any closed subscheme $\mY$ of $\mX$ and with the product structure given by
\begin{equation}
[\mX]_{k}\cdot[\mY]_{k}= [\mX\times_{k}\mY]_{k}.
\end{equation}
 
If it is clear what $k$ is, we will not write the index.
 
Denote by ${\bL} = [\mathbb{A}^1_{k}]$ the class of the affine line over $k$. By the first defining property, one obtains typical decompositions such as
\begin{equation}
[\mathbb{P}^n_k]_k = \bL^n + \cdots + \bL + 1.
\end{equation}

Let $\Sch_{\F_{q}}$ be the category of schemes of finite type over the finite field $\F_{q}$. 
An {\em additive invariant} (or ``multiplicative Euler-Poincar\'{e} characteristic'') from $\Sch_{\F_{q}}$ to $\Z$ is a map $\alpha$ which satisfies the following properties: 
\begin{itemize}
\item[(AI1)]
if $X, Y$ are in $\Sch_{\F_{q}}$ and $X \cong Y$, then $\alpha(X) = \alpha(Y)$;
\item[(AI2)]
if $X$ is in $\Sch_{\F_{q}}$ and $C$ is a closed subscheme, then $\alpha(X) = \alpha(C) + \alpha(X \setminus C)$;
\item[(AI3)]
if $X, Y$ are in $\Sch_{\F_{q}}$, then $\alpha(X \times Y) = \alpha(X)\cdot\alpha(Y)$.
\end{itemize}
An additive invariant $\alpha$ factors uniquely through a ring morphism from $K_0(\Sch_{\F_{q}})$ to $\Z$ by letting
\begin{equation}
\alpha(\sum_i[X_i]) = \sum_i\alpha[X_i].
\end{equation}

An element $\gamma$ in $K_0(\Sch_{\F_q})$, $\F_q$ a finite field, is said to be {\em polynomial-count} (according to \cite[Appendix]{Katz})  if there exists a polynomial $P(Y) = \sum_ia_iY^i \in \C[Y]$ such that for every finite extension $\F_{q^n}\Big/\F_q$, we have that $\vert \gamma(\F_q) \vert_{\F_{q^n}} = P(q^n)$. By \cite[Appendix]{Katz}, it follows that then $P(Y) \in \Z[Y]$. Note that $\vert \cdot \vert_{\F_{q^n}}$ is well defined on $K_0(\Sch_{\F_{q^n}})$ by extension from the function $\vert \cdot \vert_{\F_{q^n}}$ on $\Sch_{\F_{q^n}}$, using the fact that the latter function is an additive invariant from $\Sch_{\F_{q^n}}$ to $\Z$. 
So in terms of \cite[Appendix]{Katz}, by Theorem \ref{eqext} a $\Z$-scheme $\mX$ is of $\F_{1^\ell}$-type precisely if $[\mX_{\F_{p^\ell}}] \in K_0(\Sch_{\F_{p^\ell}})$ is polynomial-count for all primes $p$, where the corresponding counting polynomials {\em are the same for all primes $p$}. Also in terms of \cite[Appendix]{Katz}, this means that $[\mX_{p^\ell}]$ is {\em zeta equivalent} to an element $P_{\mX}$ in $\Z[\bL]$ for all primes $p$, where $P_{\mX}$ is independent of the prime $p$.

\medskip
\section{Examples}

For affine and projective spaces, we obtain the following zeta functions (over $\F_p^\ell$ and $\F_{1^\ell}$, with $\ell \in \mathbb{N}^\times$):
\begin{equation}
\left\{\begin{array}{ccc}
\zeta_{\A^n\vert \F_{p^\ell}}(s) &= &\Big({1 - p^{\ell(n - s)}}\Big)^{-1};\\
& & \\
\zeta_{\A^n\vert \F_{1^\ell}}(s) &= &\Big({\ell(s - n)}\Big)^{-1},
\end{array}\right.
\end{equation}
and
\begin{equation}\left\{
\begin{array}{ccc}
\zeta_{\bP^n\vert \F_{p^\ell}}(s) &= &\Big({(1 - p^{-\ell s})(1 - p^{\ell(1 - s)})\cdots(1 - p^{\ell(n - s)})}\Big)^{-1}; \\
& & \\
\zeta_{\bP^n\vert \F_{1^\ell}}(s) &= &\Big({\ell^{n + 1}s(s - 1)\cdots(s - n)}\Big)^{-1}.
\end{array}\right.
\end{equation}

\bigskip
We provide some more examples.
In \cite{MMKT,MMKT2}, the authors have defined and studied a functor $\mF_k$ for each field finite $k$, including $\Fun$,  which maps any ``loose graph'' $\Gamma$ (which is a more general version of a graph, where edges with $0$ or $1$ vertices are allowed) to a   $k$-scheme $\mF(\Gamma) \otimes_{\Fun}k$. Here, $\Fun$-schemes are Deitmar schemes with the additional minimal congruences considered in this paper. It is shown in \cite{MMKT} that any such scheme comes with a counting polynomial, which is in fact independent of the field. So these schemes are of $\Fun$-type, so also of $\F_{1^\ell}$-type for any positive integer $\ell$. For loose trees, precise calculations are made in \cite{MMKT} to obtain the counting polynomials. The outcome is as follows.

Let $\Gamma$ be a loose tree. We use the following notation:

\begin{itemize}
\item
$D$ is  the set of degrees $\{d_1, \ldots, d_m \}$ of $V(\Gamma)$ such that $1 < d_1 < d_2 < \ldots < d_m$;
\item
$n_i$ is the number of vertices of $\Gamma$ with degree $d_i$, $1\leq i \leq m$;
\item
$\displaystyle I= \sum_{i=1}^m n_i - 1$;
\item
$E$ is the number of vertices of $\Gamma$ with degree $1$.
\end{itemize}

Then, 
\begin{equation}
\big[\Gamma\big]_{_{\Fun}} =  \displaystyle\sum_{i = 1}^m n_i\underline{\bL}^{d_i} - I\cdot\underline{\bL} + I + E.
\end{equation}

This formula is expressed in terms of the Grothendieck ring of $\Fun$-schemes of finite type (see \cite{MMKT}). The symbol $\underline{\bL}$ denotes the class of the affine line over $\Fun$. (The same polynomial arises after base extension to fields.)

For {\em each} positive integer $\ell$, the $\F_{1^\ell}$-zeta function is thus given by
\begin{equation}
\zeta^{\F_{1^\ell}}_{\Gamma}(s)\ \ =\ \ \frac{(\ell(s - 1))^I}{(\ell s)^{E + I}}\cdot\displaystyle \prod_{k = 1}^m(\ell(s - k))^{-n_k}.
\end{equation}

If $D = \{d_1 = d \}$, $d > 1$ and $n_1 = 1$, then $I = 0$ and $E = 0$. This is the loose graph of affine $d$-space, and we obtain the $\F_{1^\ell}$-zeta function for $\A^d$ as above. (For the affine line, put $D = \emptyset$, $I = -1$ and $E = 1$.)


\newpage
\begin{thebibliography}{99}
\bibitem{ConCon}
A.~Connes and C.~Consani. Schemes over $\Fun$ and zeta functions, {\em Compos. Math.} {\bf 146} (2010), 1383--1415.


\bibitem{Deitmarschemes2}
{A.~Deitmar}.
Schemes over $\mathbb{F}_1$, In {\em Number Fields and Function Fields | Two Parallel Worlds}, Progr. Math. {\bf 239}, 2005, Birkh\"{a}user Boston, Inc., Boston, MA, pp. 87--100. 

\bibitem{Deitmarzeta}
{A.~Deitmar}.
Remarks on zeta functions and K-theory over $\mathbb{F}_1$, {\em Proc. Japan Acad. Ser. A Math. Sci.} {\bf 82} (2006), 141--146.

\bibitem{Deitmartoric}
{A.~Deitmar}.
$\mathbb{F}_1$-Schemes and toric varieties,
{\em Beitr\"{a}ge Algebra Geom.} {\bf 49} (2008), 517--525. 

\bibitem{congruence}
A.~Deitmar. Congruence schemes, {\em Internat. J. Math.} {\bf 24} (2013), 1350009, 46 pp. 


\bibitem{Deninger1991}
{C.~Deninger.}
\newblock On the {$\Gamma$}-factors attached to motives,
\newblock {\em Invent. Math.}  {\bf 104} (1991),  245--261.

\bibitem{Deninger1992}
{C.~Deninger}
\newblock Local {$L$}-factors of motives and regularized determinants,
\newblock {\em Invent. Math.} {\bf 107} (1992), 135--150.

\bibitem{Deninger1994}
{C.~Deninger.}
\newblock Motivic {$L$}-functions and regularized determinants,
\newblock In {\em Motives ({S}eattle, {WA}, 1991)}, Proc. Sympos. Pure Math. {\bf 55}, 1994, pp. 707--743.



\bibitem{Haran}
S.~Haran. Index theory, potential theory, and the Riemann hypothesis, In {\em L-functions and arithmetic (Durham, 1989)}, pp. 257--270, London Math. Soc. Lecture Note Ser. {\bf 153}, Cambridge Univ. Press, Cambridge, 1991.

\bibitem{Katz}
T.~Hausel and F.~Rodriguez-Villegas, Fernando.
Mixed Hodge polynomials of character varieties,
With an appendix by Nicholas M. Katz.
{\em Invent. Math.} {\bf 174} (2008), 555--624. 

\bibitem{KapranovUN}
{M.~Kapranov and A.~Smirnov.}
\newblock Cohomology determinants and reciprocity laws: number field case,
\newblock  Unpublished manuscript.



\bibitem{Kurokawa1992}
{N.~Kurokawa}.
\newblock Multiple zeta functions: an example,
\newblock In {\em Zeta Functions in Geometry ({T}okyo, 1990)}, Adv. Stud. Pure Math. {\bf 21}, 1992, pp. 219--226.


\bibitem{Kurozeta}
 {N.~Kurokawa}.
 Zeta functions over $\mathbf{F}_1$, {\em Proc. Japan Acad. Ser. A Math. Sci.} {\bf 81} (2005), 180--184.

\bibitem{LPL}
J.~L\'{o}pez Pe\~{n}a and O.~Lorscheid. Projective geometry for blueprints, {\em C. R. Math. Acad. Sci. Paris} {\bf 350} (2012), 455--458. 

\bibitem{LorZ}
O.~Lorscheid. Blueprints|towards absolute arithmetic? {\em J. Number Theory} {\bf 144} (2014), 408--421. 

\bibitem{Manin}
{Yu.~Manin}.
Lectures on zeta functions and motives (according to Deninger and Kurokawa), Columbia University Number Theory Seminar (New York, 1992), {\em Ast\'{e}risque} {\bf 228} (1995), 121--163.

\bibitem{MM}
Yu.~Manin and M.~Marcolli. Moduli Operad over $\Fun$, In {\em Absolute Arithmetic and $\Fun$-Geometry}, EMS Publishing House, 
30 pp., To appear.


\bibitem{MMKT}
M.~M\'{e}rida-Angulo and K.~Thas. Deitmar schemes, graphs and zeta functions, 49pp., Submitted.

\bibitem{MMKT2}
M.~M\'{e}rida-Angulo and K.~Thas. Automorphisms of Deitmar schemes, 28pp., Submitted. 

\bibitem{Soule}
{C.~Soul\'{e}}.
Les vari\'{e}t\'{e}s sur le corps \`{a} un \'{e}l\'{e}ment, {\em Mosc. Math. J.} {\bf 4} (2004), 217--244, 312.

\bibitem{KT-Chap-Comb}
K.~Thas. The Weyl functor | Introduction to Absolute Arithmetic, In {\em Absolute Arithmetic and $\Fun$-Geometry}, EMS Publishing House, 
38 pp., To appear.

\bibitem{KT-Chap-Mot}
K.~Thas. The combinatorial-motivic nature of $\Fun$-schemes, In {\em Absolute Arithmetic and $\Fun$-Geometry}, EMS Publishing House, 
84 pp., To appear.


\bibitem{KT-Chap-Weil}
K.~Thas. A taste of Weil theory, In {\em Absolute Arithmetic and $\Fun$-Geometry}, EMS Publishing House, 
28 pp., To appear.


\bibitem{anal}
{J.~Tits}.
Sur les analogues alg\'{e}briques des groupes semi-simples complexes,
In {\em Centre Belge Rech. Math., Colloque d'Alg\`{e}bre Sup\'{e}rieure, Bruxelles du 19 au 22 d\'{e}c. 1956}. 1957, pp. 261--289.






\end{thebibliography}
\end{document}